\documentclass[a4paper, english, 11pt, twoside]{amsart}

\usepackage{babel}
\usepackage[T1]{fontenc}
\usepackage{amssymb}
\usepackage{amsmath, amsfonts, mathrsfs, amscd}
\usepackage{xcolor}
\usepackage{eucal}
\usepackage{calligra}

\newtheorem{theorem}{Theorem}[section]
\newtheorem{lemma}[theorem]{Lemma}
\newtheorem{proposition}[theorem]{Proposition}
\newtheorem{corollary}[theorem]{Corollary}

\theoremstyle{definition}
\newtheorem{definition}[theorem]{Definition}
\newtheorem{remark}[theorem]{Remark}

\def\pf{\begin{proof}}
\def\epf{\end{proof}}

\newcommand{\Na}{\mathbb{N}}
\newcommand{\Z}{\mathbb{Z}}
\newcommand{\Q}{\mathbb{Q}}
\newcommand{\Co}{\mathbb{C}}
\newcommand{\Oint}{\mathcal{O}}
\newcommand{\B}{\mathcal{B}}
\newcommand{\Ind}{\operatorname{Ind}}
\newcommand{\car}{\operatorname{char}}
\newcommand{\Hom}{\operatorname{Hom}}
\newcommand{\Ch}{\mathrm{Ch}}
\newcommand{\M}{\mathrm{M}}

\newcommand{\es}{\hspace{-1pt}}

\allowdisplaybreaks

\definecolor{rojo}{rgb}{1,0,0}

\begin{document}

\title[Orders in Semisimple Hopf Algebras]{On the existence of orders in semisimple Hopf algebras}

\author[J. Cuadra and E. Meir]{Juan Cuadra and Ehud Meir}

\address{J. Cuadra: Universidad de Almer\'\i a, Dpto. Matem\'aticas, E-04120 Almer\'\i a, Spain}
\email{jcdiaz@ual.es}

\address{E. Meir: Department of Mathematical Sciences, University of Copenhagen, Universitetsparken 5, DK-2100, Denmark}
\email{meirehud@gmail.com}

\begin{abstract}
We show that there is a family of complex semisimple Hopf algebras that do not admit a Hopf order over any number ring. They are Drinfel'd twists of certain group algebras. The twist contains a scalar fraction which makes impossible the definability of such Hopf algebras over number rings. We also prove  that a complex semisimple Hopf algebra satisfies Kaplansky's sixth conjecture if and only if it admits a weak order, in the sense of Rumynin and Lorenz, over the integers.
\end{abstract}
\maketitle

\section*{Introduction}

One of the fundamental results in the Representation Theory of Finite Groups is Frobenius Theorem stating that the degree of any complex irreducible representation of a finite group $G$ divides the order of $G$, \cite[Proposition 9.32]{CR}. All known proofs of this beautiful result use a specific property of the group algebra $\Co G$:  it may be defined over the integers. In other words, $\Z G$ is an algebra order of $\Co G$. The definition of $\Co G$ from $G$ gives automatically the order $\Z G$ and, furthermore, if $\Co G$ is considered with its usual Hopf algebra structure, $\Z G$ becomes, also automatically, a Hopf order of $\Co G$. \par \smallskip

Kaplansky's sixth conjecture for Hopf algebras predicts that Frobenius Theorem holds for complex (finite dimensional) semisimple Hopf algebras, see \cite[Section 5]{So}. There are several partial results in the affirmative, listed in loc. cit. and \cite[Subsection 4.2]{A}. Compared to the group case, the main difficulty
to prove such a conjecture is that it is not guaranteed that a complex semisimple Hopf algebra $H$ is defined over the integers or, more generally, over the ring of integers of a number field. One can only ensure that $H$ is defined over a number field. Indeed, Larson proved in \cite[Proposition 4.2]{L2}, see also \cite[Theorem 4.6]{AN}, that if $H$ admits a Hopf order over a number ring, then $H$ satisfies Kaplansky's sixth conjecture (a proof is included here). The  question whether a complex semisimple Hopf algebra can be defined over a number ring was natural and has been part of the folklore of Hopf Algebra Theory. \par\smallskip

The goal of this paper is to show that, dissimilar to group algebras, a complex semisimple Hopf algebra does not necessarily admit a Hopf order over a number ring. The family of examples that we handle are Drinfel'd twists of certain group algebras. It was constructed by Galindo and Natale in \cite{GN} to produce examples of simple semisimple Hopf algebras. Let $p$ and $q$ be prime numbers such that $p$ divides $q-1$. Take an element $r \in \Z_q^{\times}$ of order $p$ and consider the semidirect product $G=(\Z_q \times \Z_q) \rtimes (\Z_p \times \Z_p),$ which can be described by generators and relations as follows:
$$\begin{array}{ll}
G= \big\langle\sigma,\tau,a,b : & \hspace{-5pt} \sigma^p=\tau^p=1,\ a^q=b^q=1,\ \sigma a\sigma^{-1}=a^r, \ \tau b \tau^{-1}=b^r , \vspace{2pt} \\
 & \hspace{-5pt} [\sigma,b]=[\tau,a]=[a,b]=[\sigma,\tau]=1\big\rangle.
\end{array}$$
Denote by $M$ the subgroup generated by $\sigma$ and $\tau$. Let $K$ be a number field containing a primitive $p$-th root of unity $\zeta$. The element
$$J=\frac{1}{p}\sum_{u,v=0}^{p-1} \zeta^{-uv}\tau^u\otimes\sigma^v \in K G \otimes K G$$ stems from a $2$-cocycle on the character group $\widehat{M}$ and it is a twist for $KG$. The twisted Hopf algebra is denoted by $\B_{p,q}(\zeta).$ In Theorem \ref{main} we prove that if $\B_{p,q}(\zeta)$ admits a Hopf order over a Dedekind domain $R \subset K$ such that $\Oint_K \subseteq R$, then $1/p \in R$. As a consequence, the complex semisimple Hopf algebra $\B_{p,q}(\zeta) \otimes_{\Q(\zeta)} \Co$ does not admit a Hopf order over any number ring (see Corollary \ref{complexif}). This result shows that the twisting operation does not preserve the definability of Hopf algebras over number rings. \par \smallskip

The proof of Theorem \ref{main} uses the following ideas and results. The elements coming from the cocharacters of a Hopf algebra must be in any Hopf order (Propositions \ref{character} and \ref{fundamental}). In this way we construct for any Hopf order $X$ of $K(\Z_q \rtimes \Z_p)$\vspace{-1.5pt} a complete set of orthogonal idempotents $\{t_i\}_{i=0}^{p-1}$ in the dual order $X^{\star}$ (Proposition \ref{classemi}). The subalgebra $A$ of $\B_{p,q}(\zeta)$ generated by $\sigma$ and $b$ is a Hopf subalgebra and $A^*\simeq K(\Z_q \rtimes \Z_p)$. Also, there is a surjective Hopf algebra map $\Phi:\B_{p,q}(\zeta)\rightarrow B$ with $B \simeq K(\Z_q \rtimes \Z_p)$. If $X$ is a Hopf order of $\B_{p,q}(\zeta)$, then the idempotents $w_0 \in X$ and $t_0 \in \Phi(X)^{\star}$, obtained from the Hopf orders $X \cap A$ of $A$ and $\Phi(X)$ of $B$ (Lemma \ref{dualorder} and Proposition \ref{subsquo}), satisfy $t_0(\Phi(w_0))=1/p$. \par \smallskip

The previous examples satisfy Kaplansky's sixth conjecture though they do not admit a Hopf order over a number ring. This reveals that, unlike group algebras, Kaplansky's sixth conjecture can not be proved by the argument that semisimple Hopf algebras are necessarily defined over a number ring. One can ask, however, if every complex semisimple Hopf algebra is isomorphic to a twist of a Hopf algebra which admits a Hopf order over a number ring. This question becomes stronger by allowing the twist operation both at the algebra level and at the coalgebra level. \par \smallskip

The essential point in proving that a complex semisimple Hopf algebra $H$ admitting a Hopf order $X$ over a number ring $R$ satisfies Kaplansky's sixth conjecture is in the Casimir element $C:=\Lambda_{(1)} \otimes S(\Lambda_{(2)})$, see Proposition \ref{kap}. Here $\Lambda$ is an integral of $H$ with $\varepsilon(\Lambda)=\dim H$. The existence of $X$ ensures that $C \in X \otimes_R X \subset H \otimes H$. Rumynin \cite{Ru} and Lorenz \cite{Lo} observed that the hypothesis on $X$ can be weakened. It is sufficient to require that $X$ is an algebra order of $H$ over $R$ and $C \in X \otimes_R X$. They call $X$ a weak order of $H$ over $R$. The same question as for Hopf orders arises: does a complex semisimple Hopf algebra admit a weak order over a number ring? \par \smallskip

In Theorem \ref{main2} we show that the answer in this context is totally different. We prove there that a complex semisimple Hopf algebra $H$ admits a weak order over $\Z$ if and only if $H$ satisfies Kaplansky's sixth conjecture. This in turn is equivalent to $C$ satisfying a monic polynomial with integer coefficients. The proof of this result relies on the decomposition of the Casimir element established in Proposition \ref{propdec}. If the Wedderburn decomposition of $H$ is $\M_{n_1}(\Co) \times \ldots \times \M_{n_s}(\Co)$, then
$$C=\sum_{i=1}^s \frac{\dim H}{n_i} C_i,$$
where $C_i$ is the Casimir element of $\M_{n_i}(\Co)$ corresponding to the bilinear form given by the trace.
\par \smallskip

The paper is organized as follows. Section 1 contains preparatory material on the Drinfel'd twist and Hopf orders. In Section 2 we discuss the Hopf orders of the group algebra of $\Z_q \rtimes \Z_p$, we define the Hopf algebra $\B_{p,q}(\zeta)$ and we prove the non existence of Hopf orders over number rings for it. Section 3 deals with weak orders. Here we first review some basic facts on Frobenius algebras, we obtain the above-mentioned decomposition of the Casimir element and we prove the announced result.
\par \smallskip

\section{Preliminaries}\label{prelim}

Throughout this paper $H$ is a finite dimensional Hopf algebra over a base field $K$. The identity element is denoted by $1_H$ and, as usual, $\Delta, \varepsilon,$ and $S$ stands for the comultiplication, counit, and antipode of $H$ respectively. For basic concepts and results on Hopf algebras we refer the reader to \cite{Mo} and \cite{Ra}. \par \smallskip

\subsection{Drinfel'd twist}
We recall from \cite[Section 2]{AEGN} and \cite[Subsection 2.4]{N} what we need about this operation.
An invertible element $J:=\sum J^{(1)} \otimes J^{(2)} \in H \otimes H$ is a {\it twist} for $H$ provided that:
$$\begin{array}{ll}
(1_H \otimes J)(id \otimes \Delta)(J)=(J \otimes 1_H)(\Delta \otimes id)(J), & \hspace{1.2cm} \textrm{{\it (2-pseudo-cocycle property)}}  \vspace{3pt} \\
(\varepsilon \otimes id)(J)=(id  \otimes \varepsilon)(J)=1_H.  & \hspace{1.2cm} \textrm{{\it (counital property)}}
\end{array}$$
Given a twist $J$ for $H$ a new Hopf algebra $H_J$, called {\it Drinfel'd twist} of $H$, can be constructed as follows: $H_J=H$ as an algebra, the counit is that of $H$, and the new comultiplication and antipode are:
\begin{equation}\label{twcomant}
\Delta_J(h)=J \Delta(h)J^{-1} \qquad \textrm{and} \qquad S_J(h)=U_JS(h)U_J^{-1} \qquad \forall h \in H.
\end{equation}
Here $U_J:=\sum J^{(1)}S(J^{(2)})$ and, if we write $J^{-1}=\sum J^{-(1)} \otimes J^{-(2)}$, its\vspace{1pt} inverse is $U_J^{-1}=\sum S(J^{-(1)})J^{-(2)}.$ \par \smallskip

We next describe the particular example of twist for the group algebra that we will use. This construction was used several times in the past, see\vspace{-1pt} for example, \cite{EG}, \cite{GN} or \cite{N} and the references therein. Let $M$ be a finite abelian group and $\widehat{M}$ its group of \hspace{0.9pt}characters. Consider \vspace{-1pt}\hspace{0.4pt} the\hspace{1.5pt} group algebra $KM$ and its dual Hopf algebra $(KM)^*$. We identify $(KM)^*$ \hspace{-0.9pt}and \hspace{-0.9pt}$K\widehat{M}$ as Hopf algebras. Assume that $\car K \nmid \vert M \vert$. There is\hspace{0.7pt} an\hspace{0.7pt} isomorphism\vspace{-1pt} of Hopf algebras $(KM)^* \simeq KM$ induced by an isomorphism of\hspace{-0.5pt} groups\hspace{-0.6pt} $\widehat{M} \simeq M$. Let $\{\vartheta_{m}\}_{m \in M} \subset (KM)^*$ be the dual basis of $\{m\}_{m \in M}$. If $\omega: M \times M \rightarrow K^{\times}$ is a normalized 2-cocycle, then $$J=\sum_{m,m' \in M} \omega(m,m') \vartheta_{m} \otimes \vartheta_{m'} \in (KM)^*\otimes (KM)^*$$ is a twist. Applying the previous isomorphism we obtain a twist for $KM$. When $G$ is an arbitrary finite group and $M$ an abelian subgroup we can construct a twist for $KG$ in this way.

\subsection{Hopf orders}
Let $R$ be a subring of $K$ and $V$ a finite dimensional $K$-vector space. An order of $V$ over $R$ is a finitely generated and projective $R$-submodule $X$ of $V$ such that the natural map $X \otimes_R K \rightarrow V$ is an isomorphism. In this case, $X$ is viewed inside $V$ as the image of $X \otimes 1_K$. A \textit{Hopf order of $H$ over $R$} is an order $X$ of $H$ such that $1_H \in X$, $XX \subseteq X$, $\Delta(X)\subseteq X\otimes_{R} X$, $\varepsilon(X) \subseteq R$ and $S(X)\subseteq X$. (Notice that $X\otimes_{R} X$ can be identified naturally as an $R$-submodule of $H\otimes H$. We shall repeatedly use this identification in the sequel.) In other words, a Hopf order of $H$ over $R$ is a Hopf algebra $X$ over $R$, which is finitely generated and projective as an $R$-module, such that $X\otimes_{R} K \simeq H$ as Hopf algebras over $K$. Observe that when $R$ is a Dedekind domain the projectivity hypothesis can be dropped in the definition: $X$ is an $R$-submodule of $H$, then it is torsion-free and hence projective. In this and the next section, and unless otherwise stated, $K$ is a number field and $R\subset K$ is a Dedekind domain which contains the ring of algebraic integers $\Oint_K$ of $K$. A Hopf order without specification of the base ring means a Hopf order over $R$. \par \smallskip

The main tool we will employ in studying Hopf orders will be that part of them is already given from the characters of the Hopf algebra and its dual. To make this idea precise, we need the dual Hopf order.
Let $X$ be a Hopf order of $H$. Recall that the dual Hopf algebra of $X$ is the dual $R$-module $X^*:=\Hom_R(X,R)$, which is finitely generated and projective, equipped with the Hopf algebra structure obtained by dualizing that of $X$. The \textit{dual order} of $X$ is defined as
$$X^{\star}=\{\varphi \in H^* : \varphi(X) \subseteq R\}.\vspace{3pt}$$
Consider the natural map $\theta:X^* \otimes_{R} K \rightarrow \Hom_K(X \otimes_{R} K,K) \simeq H^*$. Since $X$ is finitely generated and projective, $\theta$ is an isomorphism of $K$-vector spaces. One can check that $\theta(X^*  \otimes 1_K)=X^{\star}$ and $\theta$ is an isomorphism of Hopf algebras. Then:

\begin{lemma}\label{dualorder}\
\begin{enumerate}
\item[(i)] The dual order $X^{\star}$ is a Hopf order of $H^*$. \vspace{2pt}
\item[(ii)] The natural isomorphism $H \simeq H^{**}$ induces an isomorphism of Hopf algebras $X \simeq X^{\star \star}$. \qed
\end{enumerate}
\end{lemma}

We can define similarly the dual $Y^{\star}$ of an $R$-submodule $Y$ of $H$. \par \smallskip

The next result is crucial for our purposes as it states that some elements of the Hopf algebra must be in \textit{any} Hopf order.

\begin{proposition}\label{character}
Any character $\psi$ of $H$ belongs to $X^{\star}$. Similarly, any character of $H^*$ belongs to $X$. In particular, if $H$ is semisimple and $\Lambda$ is an integral with $\varepsilon(\Lambda)=\dim H$, then $\Lambda \in X$.
\end{proposition}

\pf
Suppose that $\psi$ corresponds to a representation $V$. We show that $\psi(x)\in R$ for every $x\in X$. Let $f$ be the endomorphism of $V$ induced by the action of $x$. Its eigenvalues lie in a finite extension $L$ of $K$. Since $x$ satisfies a monic polynomial with coefficients in $R$, they belong to the integral closure $\overline{R}$ of $R$ in $L$. Then the characteristic polynomial of $f$ has coefficients in $\overline{R} \cap K$. But this equals $R$ as $K$ is the field of fractions of $R$ and $R$ is integrally closed. This yields that $\psi(x) \in R$. The claim about the integral is straightforward because in a semisimple Hopf algebra the chosen integral is the character of the regular representation of $H^*$, see \cite[Theorem 16.1.2]{Ra}.
\epf

The following proposition is a consequence of \cite[Proposition 4.2]{L2}. It is also \cite[Theorem 4.6]{AN}. To keep our exposition self-contained, we include a direct proof, which follows closely the original proof of Frobenius for group algebras.

\begin{proposition}\label{kap}
Let $H$ be a split semisimple Hopf algebra over $K$. Assume that $H$ admits a Hopf order $X$ over $\Oint_K$. Then $H$ satisfies Kaplansky's sixth conjecture, i.e., the dimension of any irreducible representation of $H$ divides $\dim H$.
\end{proposition}

\pf
Let $V$ be an irreducible representation of $H$ with character $\psi$. We know that $\psi \in X^{\star}$. Let $\Lambda \in H$ be as in the previous proposition. Then $\Lambda \in X$ and, since $X$ is a Hopf order, $E:=\Lambda_{(1)}\psi(S(\Lambda_{(2)})) \in X$. We claim that $E$ is central. \vspace{-1pt} In the verification we will use that $\Lambda$ is both a left and a right integral, $S^2=id$ (by \cite[Theorem 16.1.2]{Ra}), and that characters are cocommutative. For $h\in H$ arbitrary we compute (details can be found in Lemma \ref{casimir}):
$$hE\es =\es h\Lambda_{(1)}\psi(\es S(\Lambda_{(2)})\es) \es =\es \Lambda_{(1)}\psi(S(\Lambda_{(2)})h) \es =\es \Lambda_{(1)} \psi(\es hS(\Lambda_{(2)})\es) \es =\es \Lambda_{(1)}h \psi(\es S(\Lambda_{(2)})\es) \es =\es Eh.$$
 Then $E$ must act on $V$ by a scalar $\alpha$, which belongs to $\Oint_K$ because $E \in X$. On the other hand,
$$\alpha=\frac{\psi(E)}{\dim V} = \frac{\psi (\Lambda_{(1)})\psi (S(\Lambda_{(2)}))}{\dim V} = \frac{(\psi S(\psi))(\Lambda)}{\dim V}.$$
The element $\Lambda$ equals $(\dim H)e_1$, where $e_1$ is the idempotent corresponding to the trivial representation $K$. Since the multiplicity of $K$ in $V\otimes V^*$ is $1$ ($V$ is irreducible), we have $(\psi S(\psi))(\Lambda)=\dim H$. Then $\alpha=\dim H/\dim V \in \Q \cap \Oint_K$, from which we obtain that $\alpha \in \Z$, as desired.
\epf

\begin{remark}
The proof indeed establishes that $\dim V$ divides $\dim H$ if and only if $E$ acts on $V$ by an algebraic integer. We will go deeper into this point in Subsection \ref{weak}.
\end{remark}

From the characters and cocharacters of $H$ we can form other elements as follows. For any $m,n \in \Na$ we can get linear maps $H^{\otimes m} \rightarrow H^{\otimes n}$ by using the (co)multiplication, the antipode, by applying characters of $H$ and by rearranging of factors. Evaluating such maps at $h_1\otimes\cdots\otimes h_m$, where the $h_i$'s are characters of $H^*$, we obtain elements in $H^{\otimes n}$. When $n=1$, we call them \textit{character supported elements}. For example, if $g$ and $h$ are characters of $H^*$, and $\psi$ is a character of $H$, the element $\psi(g_{(1)}g_{(4)})g_{(3)}hS(g_{(2)})$ is character supported.

\begin{definition}
Let $H$ be a split semisimple Hopf algebra over $K$. The {\it character support algebra} $\Ch_{R}(H)$ is the $R$-submodule of $H$ spanned by all the character supported elements.
\end{definition}

If $H^*$ splits over $K$, we also have the algebra $\Ch_R(H^*)$. Without loss of generality, we can assume henceforth that both $H$ and $H^*$ split over $K$ by passing, if necessary, to a finite extension of $K$, see \cite[Proposition 7.13]{CR}.

\begin{proposition}\label{fundamental}
Let $X$ be a Hopf order of $H$. Then $\Ch_R(H)\subseteq X$.
\end{proposition}

\pf
Take into account that $X$ contains all characters of $H^*$ by Proposition \ref{character}, $X$ is closed under (co)multiplication and the antipode, and if $\psi$ is a character of $H$, then $\psi(X) \subseteq R.$
\epf

The character support algebra of a group algebra $KG$ equals $RG$, and so it is already a Hopf order. However, this is not usually the case although it is possible that $\Ch_R(H)$ has maximal rank, which has the following interesting consequence:

\begin{proposition}\label{finiteorders}
Assume that $R=\Oint_K$. Any Hopf order $X$ of $H$ lies between $\Ch_R(H)$ and $\Ch_R(H^*)^{\star}$. Moreover, if these are finitely generated of maximal rank, then there are finitely many possible orders of $H$.
\end{proposition}

\pf The first claim follows from the previous lemma applied to $H$ and $H^*$. The dual order $X^{\star}$ must contain $\Ch_R(H^*)$. Hence $X \simeq X^{\star \star} \subseteq \Ch_R(H^*)^{\star}$. If $\Ch_R(H)$ and $\Ch_R(H^*)^{\star}$ are finitely generated as $R$-modules and have maximal rank $\dim H$, then the quotient $\Ch_R(H^*)^{\star}/\Ch_R(H)$ is torsion. Therefore, it is a finite abelian group and there are only finitely many options for $X$.
\epf

\begin{remark}
This result is also true under the hypotheses set on $R$ but the proof is much more involved.
\end{remark}

One last technical result about Hopf orders that we will use is the way they correspond to Hopf subalgebras and quotient Hopf algebras.

\begin{proposition}\label{subsquo}
Let $X$ be a Hopf order of a Hopf algebra $H$.
\begin{enumerate}
\item[(i)] If $A$ is a Hopf subalgebra of $H$, then $X\cap A$ is a Hopf order of $A$. \vspace{2pt}
\item[(ii)] If $\pi:H \rightarrow B$ is a surjective Hopf algebra map, then $\pi(X)$ is a Hopf order of $B$.
\end{enumerate}
\end{proposition}

\pf We first need the analogous result just at the level of orders. We recall it from \cite[Propositions 4.12 and 23.7]{CR}: \par \smallskip

Let $V$ be a finite dimensional $K$-vector space and $X$ an order of $V$. There is a bijection between the set of $R$-submodules $Y$ of $X$ such that $X/Y$ is torsion-free and the set of subspaces $W$ of $V$. It is given by $Y \mapsto KY$ and $W \mapsto X \cap W$. Moreover, $X/(X \cap W)$ is an order of $V/W$. \par \smallskip

(i) Set $Y=X \cap A$. It is clear that $1_H \in Y,\, YY\hspace{-1pt} \subseteq Y,\, \varepsilon(Y)\subseteq R,$\hspace{-1pt} and\hspace{-1pt} $S(Y) \subseteq Y$. The only difficulty is to show that $\Delta(Y) \subset Y \otimes_R Y$. Observe that $X \otimes_R X$ is an order of $H \otimes H$ and $Y \otimes_R Y$ is an order of $A \otimes A$. By the above result, $X/Y$ is torsion-free and hence projective. Then $X=Y \oplus Z$ for some $R$-submodule $Z$ of $X$. Thus $X\otimes_R X$ can be decomposed as $(Y \otimes_R Y) \oplus (Y \otimes_R Z)\oplus (Z \otimes_R Y) \oplus (Z \otimes_R Z)$. Since $A$ is a Hopf subalgebra, it must be $\Delta(Y) \subset Y \otimes_R Y$ because otherwise the $K$-vector subspace spanned by $\Delta(Y)$ would not lie inside $A \otimes A$. \par \smallskip

(ii) It is proved in a similar way. In this case, the part about the comultiplication is immediate.
\epf

\section{Non existence of Hopf orders}

In this section we construct the announced family of semisimple Hopf algebras, which are Drinfel'd twists of group algebras, that do not admit Hopf orders over \textit{any} number ring. The group algebras will be defined from a semidirect product of groups. We first need to study the Hopf orders of such algebras.

\subsection{Orders of certain group algebras}\label{classifygr}
Let $p$ and $q$ be natural numbers with $q$ prime and $p$ a divisor of $q-1$ (so that $\Z_q^{\times}$ has an element $r$ of order $p$). We take $\sigma$ and $a$ generators of $\Z_p$ and $\Z_q$ respectively. Consider the semidirect product $N:=\Z_q \rtimes\Z_p,$ where $\sigma \cdot a=a^r$. Then,
$$N = \big\langle \sigma,a : \sigma^p=a^q=1,\ \sigma a\sigma^{-1}=a^r \big\rangle.$$
Denote the group algebra $KN$ by $H$ and let $X$ be a Hopf order of $H$. From Proposition \ref{subsquo} we know that $Y:=X \cap  K \langle a \rangle$ is a Hopf order of $K\langle a \rangle$.

\begin{proposition}\label{classemi}
The Hopf order $X$ decomposes as $X=\oplus_{i=0}^{p-1} \sigma^{i}Y$.
\end{proposition}

\pf
We will show that $X^{\star}$ contains a complete set of orthogonal idempotents which will give the desired decomposition of $X$. Let $\psi$ be the character of $N$ associated to the induced representation $\Ind_{\langle \sigma \rangle}^{N} K.$ By a direct calculation\vspace{-1.7pt} one obtains \linebreak $\psi(1)=q, \psi(a^k)=0,$ and  $\psi(\sigma^ka^j)=1$ for $k\neq 0.$ By Proposition \ref{character},\vspace{-0.5pt} $\psi \in X^{\star}$. We define $\varphi \in H^*$ by $\varphi(g)=\prod_{l=0}^{p-1} \psi(\sigma^lg)$ for all $g \in N$. Then $\varphi(\sigma^ia^j)=q\delta_{j0}$. Using that $\sigma^l$ is a character of $H^*$ and $\psi(\sigma^lg)=\psi_{(1)}(\sigma^l)\psi_{(2)}(g)$, it is not difficult to realize that $\varphi \in \Ch_R(H^*)$, and consequently $\varphi \in X^{\star}$ in view of Proposition \ref{fundamental}. Define now $\nu \in H^*$ by $\nu(g)=\varphi(aga^{-1})$ for every $g \in N$. Since $X^{\star}$ is a Hopf order and $a \in X$ (Proposition \ref{character}), we have that $\nu=\varphi_{(1)}(a)\varphi_{(3)}(a^{-1})\varphi_{(2)} \in X^{\star}$. Consider\vspace{-1pt} the product $\mu:= \varphi \nu \in X^{\star}$. Let $\{\vartheta_{\sigma^ia^j}\}_{i,j} \subset H^*$ be the dual basis of $\{\sigma^ia^j\}_{i,j}$. One can easily check that $\mu = q^2 \vartheta_1$. \par \smallskip

On the other hand, the character of the regular representation of $N$ is $pq\vartheta_1$, so $pq\vartheta_1\in X^{\star}$. Since $\gcd(p,q)=1$, we get that $q\vartheta_1\in X^{\star}$. \vspace{-0.5pt} Hence $\phi:= \psi-q\vartheta_1\in X^{\star}$. For every $i=0,\ldots,p-1$ consider $t_i \in H^*$ given by
\begin{equation}\label{idemti}
t_i(g)=\prod_{l=1}^{p-1}\phi(\sigma^{l-i}g) \qquad \forall g \in N.
\end{equation}
It is easy to see that $t_i=\sum_{j=0}^{q-1}\vartheta_{\sigma^ia^j}$. \vspace{1pt} Therefore $\{t_i\}_i$ is a complete set of orthogonal idempotents in $X^{\star}$. To check that $t_i \in X^{\star}$ observe that $\phi_{(1)} \otimes \phi_{(2)} \in X^{\star} \otimes_R X^{\star}$ and $\phi_{(1)}(\sigma^{l-i})\phi_{(2)} \in X^{\star}$ as $\sigma^{l-i} \in X$. \par \smallskip

Finally, let $X^{\star}$ act on $X$ by $\gamma x = \gamma(x_{(1)})x_{(2)}$ for $\gamma \in X^{\star}, x \in X.$ Then $$X=\bigoplus_{l=0}^{p-1} t_l X.$$ If $x=\sum_{i,j} \alpha_{ij}\sigma^ia^j \in X$, where $\alpha_{ij} \in K$, then $t_l x= \sigma^l (\sum_{j} \alpha_{lj}a^j)$. Thus $t_lX=\sigma^lN_l$ for an $R$-submodule $N_l$ of $K\langle a\rangle$. Since $\sigma \in X$ and it is invertible, we have that
$N_l=N_0=X \cap K\langle a\rangle$ for all $l$. \epf

\begin{remark}
The assignment $X \mapsto X \cap  K \langle a \rangle$ establishes a bijection between the Hopf orders of $KN$ and the Hopf orders of $K\langle a \rangle$ closed under conjugation. If $K$ contains a primitive $q$-th root of unity, all Hopf orders of $K\Z_q$ are known. We review their description in \cite[Appendix]{CM} from results of Tate and Oort, Greither, and Larson. Using this description one can prove that all Hopf orders of $K\langle a \rangle$ are closed under conjugation.
\end{remark}

\subsection{The Hopf algebra $\B_{p,q}(\zeta)$}
The following family of Hopf algebras was introduced by Galindo and Natale in \cite[Section 4]{GN}. Let $p$ and $q$ be prime numbers such that $p\mid q-1$. Pick an element $r \in \Z_q^{\times}$ of order $p$. Consider the group
$$\begin{array}{ll}
G= \big\langle\sigma,\tau,a,b : & \hspace{-5pt} \sigma^p=\tau^p=1,\ a^q=b^q=1,\ \sigma a\sigma^{-1}=a^r, \ \tau b \tau^{-1}=b^r , \vspace{2pt} \\
 & \hspace{-5pt} [\sigma,b]=[\tau,a]=[a,b]=[\sigma,\tau]=1\big\rangle.
\end{array}$$
Notice that $G$ can be described as the semidirect product $(\Z_q\times \Z_q) \rtimes (\Z_p \times \Z_p)$. Also, $G=G_1 \times G_2$ with $G_1=\langle \sigma, a \rangle$ and $G_2=\langle \tau, b \rangle$. Each of these factors is isomorphic to $\Z_q \rtimes \Z_p.$ \par \smallskip

Assume that our number field $K$ contains a primitive $p$-th root of unity $\zeta$. The Hopf algebra $\B_{p,q}(\zeta)$ will be a twist of the group algebra $KG$. The twist will come from the subgroup $M$ of $G$ generated by $\sigma$ and $\tau$. Since $M$ is abelian, any twist\vspace{-1pt} for $KM$ arises from \vspace{-0.7pt} a $2$-cocycle on the character group $\widehat{M}$. For $i,j=0,\ldots,p-1$ denote by $e_{ij}$ the \vspace{-1.2pt} idempotent of $KM$ upon which $\sigma$ acts by $\zeta^i$ and $\tau$ by $\zeta^j$. It is given by:
\begin{equation}\label{idem}
e_{ij}=\frac{1}{p^2}\sum_{k,l} \zeta^{-(ik+jl)} \sigma^k\tau^l.
\end{equation}
(Throughout the limits in the sums are understood to run over all\vspace{-1.2pt} possible values.) The character in $\widehat{M}$ corresponding \vspace{-0.7pt} to $e_{ij}$ is denoted by $\psi_{ij}$. Consider\vspace{-1pt} the $2$-cocycle $\omega:\widehat{M} \times \widehat{M} \rightarrow K^{\times}, (\psi_{ij},\psi_{kl}) \mapsto \zeta^{jk}.$ The twist $J$ for $KM$ will be given by
$$J = \displaystyle{\sum\limits_{i,j,k,l} \omega(\psi_{ij},\psi_{kl})e_{ij}\otimes e_{kl}}.$$
Using \eqref{idem} one can easily check that:
$$J=\frac{1}{p}\sum_{u,v} \zeta^{-uv}\tau^u\otimes\sigma^v.$$
It is also easy to see that:
$$J^{-1} = \frac{1}{p}\sum_{u,v} \zeta^{uv}\tau^u\otimes\sigma^v.$$

Set $\B_{p,q}(\zeta)=(KG)_J$. Recall from \eqref{twcomant} that the comultiplication of $\B_{p,q}(\zeta)$ is $\Delta_J(g) = J(g\otimes g) J^{-1}$ for any $g\in G$. Clearly, $\sigma$ and $\tau$ remain group-like elements in $\B_{p,q}(\zeta)$ because $KM$ is commutative. Let us calculate $\Delta_J(a)$ and $\Delta_J(b)$:
$$\begin{array}{ll}
\Delta_J(a) & = \displaystyle{J (a\otimes a) J^{-1}} \vspace{3pt}\\
          & = \displaystyle{\frac{1}{p^2}\sum_{i,j,k,l}\zeta^{-ij+kl}\tau^i a\tau^k \otimes \sigma^j a\sigma^l} \vspace{3pt}\\
          & = \displaystyle{\frac{1}{p^2}\sum_{i,j,k,l} \zeta^{-ij+kl} a\tau^{i+k}\otimes a^{r^j}\sigma^{j+l}} \vspace{3pt} \\
          & = \displaystyle{\frac{1}{p^2}\sum_{u,j,l} \zeta^{-uj}\bigg(\sum_{k} \zeta^{k(j+l)}\bigg)a\tau^{u}\otimes a^{r^j}\sigma^{j+l}} \vspace{3pt} \\
          & = \displaystyle{\frac{1}{p}\sum_{u,j} \zeta^{-uj} a\tau^u\otimes a^{r^j}.}
\end{array}$$
Similarly,
\begin{equation}\label{comulty}
\Delta_{J}(b) = \frac{1}{p}\sum_{i,v} \zeta^{-iv}b^{r^i}\otimes b\sigma^v.
\end{equation}
Recall also from \eqref{twcomant} that the antipode of $\B_{p,q}(\zeta)$ is $S_J(g) = U_JS(g)U_J^{-1}$ for any $g\in G$. We compute $S_J(a)$ and $S_J(b)$:
$$\begin{array}{ll}
S_J(a)    & = \displaystyle{U_J a^{-1} U_J^{-1}} \vspace{3pt}\\
          & = \displaystyle{\frac{1}{p^2}\sum_{i,j,k,l}\zeta^{-ij+kl} \tau^i \sigma^{-j} a^{-1} \tau^{-k} \sigma^l} \vspace{3pt}\\
          & = \displaystyle{\frac{1}{p^2}\sum_{i,j,k,l} \zeta^{-ij+kl} a^{-r^{-j}}\tau^{i-k}\sigma^{l-j}} \vspace{3pt} \\
          & = \displaystyle{\frac{1}{p^2}\sum_{u,v,l} \zeta^{u(v-l)} \bigg(\sum_{k} \zeta^{vk}\bigg)a^{-r^{v-l}}\tau^{u} \sigma^{v}} \vspace{3pt} \\
          & = \displaystyle{\frac{1}{p} \sum_{u,l} \zeta^{-ul} a^{-r^{-l}}\tau^u.}
\end{array}$$
Similarly,
\begin{equation}\label{antipy}
S_J(b) = \frac{1}{p} \sum_{k,v} \zeta^{vk} b^{-r^{k}}\sigma^v.
\end{equation}

\subsection{The Hopf algebra $\B_{p,q}(\zeta)$ has no Hopf orders}
We are now in a position to state and prove the main result of this paper:

\begin{theorem}\label{main}
Let $K$ be a number field containing a primitive $p$-th root of unity $\zeta$. Let $R\subset K$ be a Dedekind domain such that $\Oint_K \subseteq R$. If $\B_{p,q}(\zeta)$ admits a Hopf order over $R$, then $\frac{1}{p}\in R$. As a consequence, $\B_{p,q}(\zeta)$ does not admit a Hopf order over any number ring.
\end{theorem}

\pf
Consider the subalgebra $A$ of $\B_{p,q}(\zeta)$ generated by $\sigma$ and $b$, which is isomorphic to $K(\Z_p \times \Z_q)$. By \eqref{comulty} and \eqref{antipy}, $A$ is a Hopf subalgebra of $\B_{p,q}(\zeta)$. Denote by $N$ the character group of $A$. We can assume, without loss of generality, that $K$ contains a primitive $q$-th root of unity $\eta$. Then $A$ splits over $K$. Consider
$$f_{kl}=\frac{1}{pq}\sum_{i}\sum_{j} \zeta^{ik}\eta^{-lj} \sigma^ib^j, \qquad k=0,\ldots, p-1, \ l=0,\ldots, q-1,$$
the complete set of orthogonal idempotents giving the Wedderburn decomposition of $A.$ Let $\psi_{kl}$ be the character of $A$ corresponding to $Kf_{kl}$. It is given by $\psi_{kl}(\sigma^ib^j)=\zeta^{-ki}\eta^{lj}$. We have that $A^*$ is isomorphic to $KN$ as a Hopf algebra. We claim that $N$ is isomorphic to $\Z_q \rtimes \Z_p$. Set $s=r^{-1}\hspace{-5pt} \mod q$. We compute:
$$\begin{array}{ll}
(\psi_{kl}\psi_{cd})(\sigma) & = \psi_{kl}(\sigma)\psi_{cd}(\sigma) \vspace{3pt} \\
                            & = \zeta^{-(k+c)} \vspace{3pt} \\
                            & = \psi_{k+c\hspace{2pt} ls^c+d}(\sigma). \vspace{9pt} \\
(\psi_{kl}\psi_{cd})(b) & = \displaystyle{\frac{1}{p}\sum_{i,v} \zeta^{-iv} \psi_{kl}(b^{r^i})\psi_{cd}(b\sigma^{v})} \vspace{3pt} \\
                        & = \displaystyle{\frac{1}{p} \sum_{i} \eta^{lr^i+d} \bigg(\sum_{v} \zeta^{-(c+i)v}\bigg)} \vspace{3pt}\\
                        & = \eta^{lr^{-c}+d} \\
                        & = \psi_{k+c\hspace{2pt}ls^c+d}(b).
\end{array}$$
Then $\psi_{10}^p=\psi_{01}^q=1, \psi_{10}\psi_{01}\psi_{10}^{-1}=\psi_{01}^r.$ \par \smallskip

Let $X$ be a Hopf order of $\B_{p,q}(\zeta)$ over $R$. By Proposition \ref{subsquo} (i) and Lemma \ref{dualorder} (i), $(X \cap A)^{\star}$ is a Hopf order of $A^* \simeq K(\Z_q \rtimes \Z_p)$. We have seen in the proof of Proposition \ref{classemi} that $X \cap A \simeq (X \cap A)^{\star \star}$ must contain the idempotents $t_i$ in \eqref{idemti} for $i=0,\ldots,p-1$.
To avoid confusion later we denote them here by $w_i$. Since ${(f_{kl},\psi_{kl})}_{k,l}$ is a dual basis of $A^*$, the idempotent $w_i$ is given by
$$w_i=\sum_{j} f_{ij} = \frac{1}{p}\sum_{u} \zeta^{iu}\sigma^u.$$
In particular, $w_0=\frac{1}{p}\sum_{u} \sigma^u \in X \cap A$. \par \smallskip

Now, since $G=G_1 \times G_2$, we have a surjective Hopf algebra map $\Phi':KG\rightarrow KG_1$ given by the projection of groups $G \twoheadrightarrow G_1, \sigma^i\tau^ja^kb^l \mapsto \sigma^ia^k$. This induces a Hopf algebra map $\Phi:\B_{p,q}(\zeta) \rightarrow (KG_1)_{\bar{J}},$ where $\bar{J}=(\Phi' \otimes \Phi')(J)$. It is easy to see that $\bar{J}=1\otimes 1$, so the twist is trivial. By Proposition \ref{subsquo} (ii), $\Phi(X)$ is a Hopf order of $KG_1$. Recall that $G_1 \simeq \Z_q \rtimes \Z_p$. This order must contain $$y:=\Phi\bigg(\frac{1}{p}\sum_{u} \sigma^u\bigg) = \frac{1}{p}\sum_{u} \sigma^u.$$
Consider $t_0 \in \Phi(X)^{\star}$ as in \eqref{idemti}. Then $t_0(y)=\frac{1}{p} \in R.$ \epf

\begin{corollary}\label{complexif}
The complex semisimple Hopf algebra $\B_{p,q}(\zeta) \otimes_{\Q(\zeta)} \Co$ does not admit a Hopf order over any number ring.
\end{corollary}

\pf
Suppose that $\B_{p,q}(\zeta) \otimes_{\Q(\zeta)} \Co$ admits a Hopf order $X$ over a number ring $\Oint_K$. We can assume that $\Q(\zeta) \subseteq K$. Then we would have
$(\B_{p,q}(\zeta) \otimes_{\Q(\zeta)} K) \otimes_K \Co \simeq (X \otimes_{\Oint_K} K) \otimes_K \Co.$
This implies that $\B_{p,q}(\zeta) \otimes_{\Q(\zeta)} K$ and $X \otimes_{\Oint_K} K$ would become isomorphic when extending them by scalars to some finite  extension $L$ of $K$, and $X \otimes_{\Oint_K} \Oint_L$ would be a Hopf order of $\B_{p,q}(\zeta) \otimes_{\Q(\zeta)} L$ over $\Oint_L$, violating Theorem \ref{main}. \epf

\begin{remark}
We can construct more examples of Hopf algebras that do not admit Hopf orders over number rings out of the previous ones by taking duals, Drinfel'd double or Hopf algebra extensions. Let $B$ and $H$ be Hopf algebras and suppose that $B$ does not admit a Hopf order over any number ring. If either $B$ is a Hopf subalgebra or a quotient of $H$, then $H$ does not admit a Hopf order over any number ring. This follows from Proposition \ref{subsquo}. For the statement on duals take into account Lemma \ref{dualorder}.
\end{remark}

\subsection{Minimal example}
Notice that the smallest example we have constructed of a Hopf algebra without a Hopf order is $\B_{2,3}(-1)$, whose dimension is $36$. We wonder if this is the minimal example. We can check directly, using known classification results, that any complex semisimple Hopf algebra of dimension less than $36$ and different from $24$ and $32$ (where the classification is open) has a Hopf order over the ring of integers of a cyclotomic field. \par \smallskip

Recall that a Hopf algebra is called trivial if it is isomorphic to either a group algebra or the dual of a group algebra. A trivial Hopf algebra admits a Hopf order over $\Z$. By well-known results of Etingof and Gelaki, Masuoka, Natale, and Zhu, a complex semisimple Hopf algebra $H$ is trivial if:
\begin{enumerate}
\item $\dim H$ is prime.
\item $\dim H=p^2$, with $p$ prime.
\item $\dim H=pq$, with $p$ and $q$ two distinct prime numbers.
\item $\dim H=30$.
\end{enumerate}

Semisimple Hopf algebras of dimension $p^3$ were classified by Masuoka in \cite{M1} and \cite{M2}. The non-trivial examples admit a Hopf order over either $\Z$ or $\Z(\zeta)$, with $\zeta$ a primitive $p$-th root of unity, see \cite[Remark 2.14 (1)]{M1} and \cite[Section 2]{M2}. \par \smallskip

Semisimple Hopf algebras of dimension $pq^2$ were classified by Natale. Since we are dealing with the case $pq^2<36$, the partial classification established in \cite{N1} suffices for us. The non-trivial examples are defined over $\Z(\zeta)$, with $\zeta$ a root of unity whose order depends of each case, see \cite[Lemma 1.3.11 and (1.4.9)]{N1} and \cite[page 249]{G}. \par \smallskip

Semisimple Hopf algebras of dimension $16$ were classified by Kashina in \cite{K}. Any such Hopf algebra $H$ is a cleft extension $A \#_{\sigma}^{\theta}\, \Co\Z_2,$ where $A$ is a normal commutative Hopf subalgebra of dimension $8$, see \cite[Theorem 9.1]{K}. Then $A \simeq (\Co G)^*$ for a group $G$ of order $8$. Let $\{\vartheta_g\}_{g \in G} \subset A$ be the dual basis and $t$ a generator of $\Z_2$. A case by case study reveals that $\{\vartheta_g \# t^j: g\in G, j=0,1\}$ is an integral basis for $H$. The structure constants in this basis belong to $\Z(\zeta)$, with $\zeta$ a primitive $8$-th root of unity. \par \smallskip

We refer the reader to \cite{N2} for more results on the classification of low dimensional semisimple Hopf algebras. It is a curious coincidence that $\B_{2,3}(-1)$ is the minimal example of a non-trivial simple semisimple Hopf algebra.

\section{Weak orders}

There is a weak notion of Hopf order, due to Rumynin \cite{Ru} and studied further by Lorenz \cite{Lo}, under which Kaplansky's sixth conjecture still holds. The same line of thought lead to ask if every complex semisimple Hopf algebra admits a weak order over a number ring. In this final section we prove that the existence of a weak order over $\Z$ for a complex semisimple Hopf algebra $H$ amounts to $H$ satisfies Kaplansky's sixth conjecture. \par \smallskip

We start by recalling the structure of Frobenius algebra of any finite dimensional Hopf algebra and the description of the Casimir element. In this section $K$ is an arbitrary field.

\subsection{Hopf algebras and Frobenius algebras}
Recall from \cite[Proposition 9.5]{CR} that a finite dimensional $K$-algebra $A$ is called {\it Frobenius} if there exists a nondegenerate associative bilinear form $\langle \cdot, \cdot \rangle: A \times A \rightarrow K.$ The {\it Frobenius homomorphism} of $A$ is the map $\phi:A \rightarrow K, a \mapsto \langle a, 1_A \rangle$. Given $a \in A$ consider the map $\phi(a\hspace{2pt} \cdot):A \rightarrow K, b \mapsto \phi(ab).$ The nondegeneracy of $\langle \cdot, \cdot \rangle$ implies that $\Phi: A \rightarrow A^*, a  \mapsto \phi(a)$ is an isomorphism of $K$-vector spaces. Let $\{a_i\}_i$ be a basis of $A$ and $\{a_i^*\}_i \subset A^*$ its dual basis. Set $b_i=\Phi^{-1}(a_i^*)$. The {\it Casimir element} of $A$ is $c:=\sum_i a_i \otimes b_i$. It is uniquely determined by the property $\sum_i \langle b_i, x \rangle a_i =x$ for all $x \in A$. Equivalently, $\sum_i \langle x, a_i \rangle b_i=x$ for all $x \in A$. It is known that $(a \otimes 1_A)c=c(1_A \otimes a)$ for all $a \in A$. \par \smallskip

The matrix algebra $\M_n(K)$ is Frobenius with bilinear form given by the trace, i.e., $\langle \cdot, \cdot \rangle: \M_n(K) \times \M_n(K) \rightarrow K, (P,Q) \mapsto tr(PQ).$ For $i,j=1,\ldots,n$ we denote by $T_{ij}$ the matrix with $1$ in the $(i,j)$-entry and zero elsewhere. The Casimir element of $\M_n(K)$ is $\sum_{k,l=1}^n T_{kl} \otimes T_{lk}$. \par \smallskip

Let now $H$ be a finite dimensional Hopf algebra over $K$. Take $\lambda \in H^*$ a nonzero right integral and $\Lambda \in H$ a left integral such that $\lambda(\Lambda)=1$. The bilinear form $\langle \cdot, \cdot \rangle: H \times H \rightarrow K, (g,h) \mapsto \lambda(gh)$ is associative and nondegenerate and the corresponding Casimir element is $C:=\Lambda_{(1)} \otimes S(\Lambda_{(2)}).$ Similarly, if $\gamma \in H^*$ is a nonzero left integral and $\Gamma \in H$ a right integral with $\gamma(\Gamma)=1$, the bilinear form $\langle \cdot, \cdot \rangle: H \times H \rightarrow K,$ \linebreak $(g,h) \mapsto \gamma(gh)$ is  associative and nondegenerate and the Casimir element is \linebreak $S(\Gamma_{(1)}) \otimes \Gamma_{(2)}.$ \par \smallskip

The following result is fundamental for our aim:

\begin{lemma}\label{casimir}
Let $H$ be a finite dimensional unimodular Hopf algebra over $K$. Then, $H$ is involutory if and only if $(g \otimes h)C=C(h \otimes g)$ for all $g,h \in H.$
\end{lemma}

\pf As we pointed out before, $(g \otimes 1_H)C=C (1_H \otimes g)$ holds in any Frobenius algebra. Assume that $H$ is involutory. That $(1_H \otimes h)C=C(h \otimes 1_H)$ follows from \cite[Proposition 10.1.3 (b)]{Ra} using that $S^2=id$. We include the proof for completeness. Since $H$ is unimodular, $\Lambda$ is also a right integral. We compute:
$$\begin{array}{ll}
C(h \otimes 1_H) & = \Lambda_{(1)} h \otimes S(\Lambda_{(2)}) \vspace{3pt} \\
               & = \Lambda_{(1)} h_{(1)} \otimes \varepsilon(h_{(2)}) S(\Lambda_{(2)}) \vspace{3pt} \\
               & = \Lambda_{(1)} h_{(1)} \otimes h_{(3)}S(h_{(2)}) S(\Lambda_{(2)}) \vspace{3pt} \\
\end{array}$$
$$\begin{array}{ll}
               & = \Lambda_{(1)} h_{(1)} \otimes h_{(3)} S(\Lambda_{(2)}h_{(2)}) \vspace{3pt} \\
               & = \Lambda_{(1)} \otimes \varepsilon(h_{(1)})h_{(2)}S(\Lambda_{(2)})  \vspace{3pt} \\
               & = \Lambda_{(1)} \otimes h S(\Lambda_{(2)}) \vspace{3pt} \\
               & = (1_H \otimes h)C .
\end{array}$$
Conversely, suppose that $(g \otimes h)C=C(h \otimes g)$ for all $g,h \in H.$ By \cite[Proposition 10.1.3 (b)]{Ra}, $\Lambda_{(1)} S^{-1}(h) \otimes \Lambda_{(2)}=\Lambda_{(1)} \otimes \Lambda_{(2)}h$. Then,
$$\Lambda_{(1)} S^{-1}(h) \otimes S(\Lambda_{(2)})=\Lambda_{(1)} \otimes S(h)S(\Lambda_{(2)})=\Lambda_{(1)} S(h) \otimes S(\Lambda_{(2)}).$$
From here,
$$hS(\Lambda_{(1)}) \otimes \Lambda_{(2)}=S^2(h) S(\Lambda_{(1)}) \otimes \Lambda_{(2)}.$$
Since $S(\Lambda_{(1)}) \otimes \Lambda_{(2)}$ is the Casimir element for the bilinear form arising from a nonzero left integral, we have:
$$h=\langle h, S(\Lambda_{(1)}) \rangle \Lambda_{(2)} = \langle S^2(h), S(\Lambda_{(1)})\rangle \Lambda_{(2)}=S^2(h).$$
\epf

\subsection{Weak orders}\label{weak} See \cite{Ru} and \cite[Subsections 2.5 and 5.9 and Corollary 16]{Lo}. Let $H$ be a finite dimensional Hopf algebra over $K$ and $R$ a subring of $K$. A {\it weak order} of $H$ over $R$ is an algebra order $X$ of $H$ such that the Casimir element $C \in H \otimes H$ lies in $X \otimes_R X$. Recall that $X\otimes_{R} X$ is viewed naturally as an $R$-submodule of $H\otimes H$. \par \smallskip

Observe that this definition captures the key point in the proof of Proposition \ref{kap}. There the Hopf order is used precisely to ensure that the Casimir element fulfill this. The proof actually establishes the following result. Let $H$ be a split semisimple Hopf algebra over a number field $K$ and assume that $H$ admits a weak order $X$ over $\Oint_K$. Then $H$ satisfies Kaplansky's sixth conjecture. Notice that in proving that $\psi(x) \in \Oint_K$ for all $x \in X$ we only need that $X$ is an algebra order, see Proposition \ref{character}.

\subsection{Decomposition of the Casimir element of a semisimple Hopf algebra}
We are going to see here that in a split, semisimple, and involutory Hopf algebra the Casimir element can be described from the Casimir elements attached to the trace of each matrix constituent. The latter are uniquely determined, up to scalar, by the following property:

\begin{proposition}\label{matrix}
An element $D \in \M_n(K) \otimes \M_n(K)$ satisfies $(P \otimes Q)D=D(Q \otimes P)$ for all $P,Q \in \M_n(K)$ if and only if
$D=\beta \sum_{k,l=1}^n T_{kl} \otimes T_{lk}$ for some $\beta \in K$. In this case, $D^2=\beta^2 (1 \otimes 1)$.
\end{proposition}

\pf Write $D=\sum_{k,l,u,v=1}^n \beta_{kl}^{uv}T_{kl} \otimes T_{uv}$ with $\beta_{kl}^{uv} \in K$. We calculate:
\begin{align*}
(T_{ij} \otimes T_{cd})D & = \sum_{k,l,u,v=1}^n \beta_{kl}^{uv}T_{ij} T_{kl} \otimes T_{cd}T_{uv}= \sum_{l,v=1}^n \beta_{jl}^{dv}T_{il} \otimes T_{cv}. \\
D(T_{cd} \otimes T_{ij}) & = \sum_{k,l,u,v=1}^n \beta_{kl}^{uv} T_{kl}T_{cd} \otimes T_{uv}T_{ij}= \sum_{k,u=1}^n \beta_{kc}^{ui}T_{kd} \otimes T_{uj}.
\end{align*}
If either $u \neq c$ or $v \neq j$, then $\beta_{kc}^{ui}=\beta_{jl}^{dv}=0$. Otherwise, $\beta_{kc}^{ci}=\beta_{jl}^{dj}$. Setting $\beta=\beta_{ic}^{ci}$, we have:
$$D=\beta\sum_{k,l=1}^n T_{kl} \otimes T_{lk}.$$
For the second statement, we calculate:
$$D^2=\beta^2 \sum_{k,l,u,v=1}^n T_{kl}T_{uv} \otimes T_{lk}T_{vu}= \beta^2\sum_{u,v=1}^n T_{vv} \otimes T_{uu}= \beta^2 (1 \otimes 1).$$
\epf

\begin{proposition}\label{propdec}
Let $H$ be a split semisimple Hopf algebra over $K$ and $0 \! \neq \! \Lambda \! \in \! H$ an integral. Assume that the Wedderburn decomposition of $H$ is $\M_{n_1}(K) \times \ldots \times \M_{n_s}(K)$. The following assertions are equivalent:
\begin{enumerate}
\item[(i)] $H$ is involutory. \vspace{5pt}
\item[(ii)] The Casimir element $C=\Lambda_{(1)} \otimes S(\Lambda_{(2)})$ can be written as
\begin{equation}\label{decomp}
C=\sum_{i=1}^s \frac{\varepsilon(\Lambda)}{n_i} C_i,
\end{equation}
where $C_i$ is the image of the Casimir element of $\M_{n_i}(K)$ associated to the trace.
\end{enumerate}
\end{proposition}

\pf Let $\{e_i\}_{i=1}^s$ be the complete set of orthogonal central idempotents attached to the Wedderburn decomposition of $H$. We know that $He_i=e_iH \simeq \M_{n_i}(K)$. Put $f_{ij}=e_i \otimes e_j$.  Then $\{f_{ij}\}_{i,j=1}^s$ is the complete set of orthogonal central idempotents attached to the Wedderburn decomposition of $H \otimes H$. \par \smallskip

(i) $\Rightarrow$ (ii) We know from Lemma \ref{casimir} that $(g \otimes h)C=C(h \otimes g)$ for all $g,h \in H$. Then
$$\begin{array}{ll}
C & = \sum_{i=1}^s C(1_H \otimes e_i) \vspace{3pt} \\
  & = \sum_{i=1}^s C(1_H \otimes e_i)(1_H \otimes e_i) \vspace{3pt} \\
  & = \sum_{i=1}^s (e_i \otimes 1_H)C(1_H \otimes e_i) \vspace{3pt} \\
  & = \sum_{i=1}^s (e_i \otimes e_i)C \vspace{3pt} \\
  & = \sum_{i=1}^s f_{ii} C.
\end{array}$$
Now, $f_{ii}C \in f_{ii}(H \otimes H)=e_iH \otimes e_i H \simeq \M_{n_i}(K) \otimes \M_{n_i}(K)$. For $g,h \in H$ arbitrary, we have
$(ge_i \otimes he_i)f_{ii}C=Cf_{ii}(he_i \otimes ge_i)$. By Proposition \ref{matrix}, $f_{ii}C=\beta_iC_i$ for some $\beta_i \in K$, where $C_i$ denotes the Casimir element of $\M_{n_i}(K)$ associated to the trace. Let $d^i_{kl} \in e_iH$ denote the element corresponding to the matrix $T_{kl}$ in $\M_{n_i}(K).$ Then  $f_{ii}C=\beta_i \sum_{k,l=1}^{n_i} d^i_{kl} \otimes d^i_{lk}$. On the other hand,
$$ \sum_{i=1}^s \varepsilon(\Lambda)e_i = \varepsilon(\Lambda)1_H = \Lambda_{(1)}S(\Lambda_{(2)}) =  \sum_{i=1}^s \beta_i \sum_{k,l=1}^{n_i} d^i_{kl}d^i_{lk} = \sum_{i=1}^s \beta_in_ie_i.$$
Thus $\beta_in_i=\varepsilon(\Lambda)$ and we obtain the desired decomposition of $C$. \par \smallskip

(ii) $\Rightarrow$ (i) By Proposition \ref{matrix}, $(ge_i \otimes he_i)C_i=C_i(he_i \otimes ge_i)$ for all $g,h \in H.$ This yields that
$(g \otimes h)C=C(h \otimes g)$. By Lemma \ref{casimir}, $H$ is involutory.
\epf

We have set up all the necessary elements to prove the main result of this section.

\begin{theorem}\label{main2}
Let $H$ be a split semisimple Hopf algebra over a field of characteristic zero. Let $\Lambda \in H$ be an integral with $\varepsilon(\Lambda)=\dim H$. The following assertions are equivalent:
\begin{enumerate}
\item[(i)] $H$ admits a weak order over $\Z.$
\item[(ii)] The Casimir element $\Lambda_{(1)} \otimes S(\Lambda_{(2)})$ satisfies a monic polynomial with coefficients in $\Z.$
\item[(iii)] $H$ satisfies Kaplansky's sixth conjecture.
\end{enumerate}
\end{theorem}

\pf By Larson-Radford's Theorem \cite[Theorem 16.1.2]{Ra}, $H$ is involutory. As before, call $C$ the Casimir element. \par \smallskip

(i) $\Rightarrow$ (ii) If $H$ admits a weak order $X$ over $\Z$, then $C \in X \otimes_{\Z} X \subset H \otimes H$. Since the $\Z$-algebra $X \otimes_{\Z} X$ is finitely generated as a $\Z$-module, the element $C$ satisfies a monic polynomial with coefficients in $\Z.$ \par \smallskip

(ii) $\Rightarrow$ (iii)  Assume that $C$ satisfies a monic polynomial $m(t) \in \Z[t].$ We keep the notation of the previous proof. Consider the projection $\pi_i: H  \otimes H \rightarrow (H \otimes H)f_{ii}$  on the Weddeburn component attached to $f_{ii}$. By \eqref{decomp}, $\pi_{i}(C)=n_i^{-1}(\dim H) C_{i}$. Hence $n_i^{-1}(\dim H) C_{i}$ also satisfies $m(t)$. Recall that $C_i^2=f_{ii} \otimes f_{ii}$ by Proposition \ref{matrix}. The above implies that $n_i^{-1}\dim H$ is an algebraic integer and, therefore, an integer. \par \smallskip

(iii) $\Rightarrow$ (i) The $\Z$-algebra $Y:= \M_{n_1}(\Z) \times \ldots \times \M_{n_s}(\Z)$ is an algebra order of $H$ over $\Z$. If $H$ satisfies Kaplansky's sixth conjecture, $n_i \mid \dim H$ for all $i$. The decomposition \eqref{decomp} of $C$ implies that $C \in Y \otimes_{\Z} Y$ and we are done.
\epf

\begin{remark}
The above proof tells us how to construct a weak order over the integers for any complex semisimple Hopf algebra satisfying Kaplansky's sixth conjecture. The weak order is given by the Wedderburn decomposition where each matrix block has only integer coefficients. In particular, it is not difficult to construct a weak order for the Hopf algebra $\B_{p,q}(\zeta)$. As an algebra, $\B_{p,q}(\zeta)$ is isomorphic to the group algebra of $(\Z_q \rtimes \Z_p) \times (\Z_q \rtimes \Z_p)$. Its Wedderburn decomposition appears for example in \cite[Lemma 4.1]{GN}:
$$\B_{p,q}(\zeta) \simeq K^{(p^2)} \times \underbrace{\M_p(K) \times \ldots \times \M_p(K)}_{2(q-1) \ \text{copies}} \times \underbrace{\M_{p^2}(K) \times \ldots \times \M_{p^2}(K)}_{\big(\frac{q-1}{p}\big)^2\ \text{copies}}.$$
The characters of $\Z_q \rtimes \Z_p$ can be calculated using \cite[Theorem 5.27.1]{EEA} and the complete set of orthogonal central idempotents using \cite[Proposition 9.21]{CR}.
\end{remark}

\vspace{7pt}

\subsection*{Acknowledgements}
The first author was supported by the projects MTM2011-27090 from MICINN and FEDER and by the research group FQM0211 from Junta de Andaluc\'{\i}a. The second author was supported by the Danish National Research Foundation (DNRF) through the Centre for Symmetry and Deformation. \par \smallskip

The authors would like to thank Dmitriy Rumynin and Martin Lorenz for bringing to our attention the notion of weak order; C\'esar Galindo for pointing out that these examples were already discussed in \cite{GN}; Yevgenia Kashina and Sonia Natale for useful discussions on the classification of low dimensional semisimple Hopf algebras.

\end{document}